\documentclass[10pt]{amsart}

\newtheorem{theorem}{Theorem}[section]
\newtheorem{lemma}[theorem]{Lemma}
\newtheorem{proposition}[theorem]{Proposition}
\newtheorem{corollary}[theorem]{Corollary}
\theoremstyle{definition}
\newtheorem{definition}{Definition}[section]
\theoremstyle{remark}
\newtheorem{remark}{Remark}[section]
\theoremstyle{remark}

\usepackage{amssymb}
\newcommand{\vertt}{|\!|\!|}
\renewcommand{\leq}{\leqslant}
\renewcommand{\geq}{\geqslant}
\renewcommand{\le}{\leqslant}
\renewcommand{\ge}{\geqslant}
\renewcommand{\L}{\mathrm{L}}
\newcommand{\C}{\mathrm{C}}
\newcommand{\e}{\,\mathrm{e}\,}
\newcommand{\ds}{\displaystyle}
\newcommand{\R}{\mathbb{R}}
\renewcommand{\d}{\,\mathrm{d}}

\begin{document}

\title[Isothermalisation for a Non-local Heat Equation]
 {Isothermalisation for a Non-local\\ Heat Equation}

\author[Chasseigne]{Emmanuel Chasseigne}
\address{Laboratoire de Math\'ematiques et Physique Th\'eorique,\\
Universite de Tours, Parc de Grandmont, F37200 Tours\\
FRANCE\\}
\email{emmanuel.chasseigne@lmpt.univ-tours.fr}

\author[Ferreira]{Ra\'ul Ferreira}
\address{Departamento de Matem\'atica Aplicada,
Fac. de C.C. Qu\'{\i}micas,\\
U.~Complutense de Madrid, 28040 Madrid\\
SPAIN}
\email{raul$_-$ferreira@mat.ucm.es}



\begin{abstract}
In this paper we study the asymptotic behavior for a nonlocal heat
equation in an inhomogenous medium:
$$\rho(x)u_t=J\ast u-u \text{ in }\mathbb{R}^N\times (0,\infty)\,,$$
where $\rho$ is a continous positive function, $u$ is nonnegative and $J$ is a probability measure
having finite second-order momentum. Depending on integrability conditions on the initial data $u_0$
and $\rho$, we prove various isothermalisation results, \textit{i.e.}
$u(t)$ converges to a constant state in the whole space.
\end{abstract}

\maketitle

\section{Introduction}

The aim of this paper is to study the asymptotic behavior for a nonlocal heat
equation in an inhomogenous medium:
\begin{equation}\label{eq:isot.0}
    \left\{
    \begin{array}{ll}
    \rho(x)u_t=J\ast u-u, \quad &  (x,t)\in \mathbb{R}^N\times (0,\infty),\\
    u(x,0)=u_0(x), &x\in\mathbb{R}^N.
    \end{array}\right.
\end{equation}
Here, $u_0$ is a nonnegative continuous function in $\R^N$ and
$\ast$ denotes the convolution with a kernel $J:\mathbb R^N \to
\mathbb R$, which is a radial, continuous probability density having
finite second-order momentum:
$$
    \int_{\mathbb{R}^N} J(s)\d s=1\,,\quad \mathbb{E}(J)=\int_{\mathbb{R}^N} s J(s)\d s=0\,,\
    \mathbb{V}(J)=\int_{\mathbb{R}^N} s^2 J(s)\d s<+\infty\,.
$$
Typical examples of kernel that we consider are the gaussian law,
the exponential law or any compactly supported kernels. We also
assume that $\rho$ is a positive, continuous function
in $\R^N$, whether integrable or not.

The operator $J*u-u$ can be interpreted as a non-local diffusion
operator. Indeed, if $u(x,t)$ represents the density of a single
population and $J(x-y)$ is the probability to jump from $y$ to $x$
then the term $(J*u)(x)$ is  the rate at which individuals arrive to
$x$ and $-u(x)$ is the rate at which individuals leave from $x$, see
for instance \cite{Fi}. In the case of heat propagation, $u$ stands
for a temperature and $\rho(x)$ represents the density of the
medium.

Problem (\ref{eq:isot.0}) is called {\em non-local} because the
diffusion at $u(x,t)$ depends on all the values of $u$ in the
support of $J$ and not only of the value of $u(x,t)$, as it is the case for
the local diffusion problem
$$
\left\{
    \begin{array}{ll}
    \rho(x)u_t=\Delta u, \quad &  (x,t)\in \mathbb{R}^N\times (0,\infty),\\
    u(x,0)=u_0(x), &x\in\mathbb{R}^N.
    \end{array}
\right.
$$
For this local problem it is well known that for dimension $N=1,2$
there exists a unique solution in the class of bounded solutions,
see \cite{KR} and \cite{Eidus}.
Moreover, if $\rho\in\L^1(\R^N)$ and $u_0$ is bounded, then as $t\to\infty$ the solution converges on compact sets to
$\mathbb{E}_\rho(u_0)$, the mean of $u_0$ with respect to $\rho$~:
$$
    \mathbb{E}_\rho(u_0):=
    \frac{\ds\int_{\R^N}u_0(y)\rho(y)\d y}{\ds\int_{\R^N}\rho(y)\d y}\in\R_+\,.
$$
This phenomenon is called {\em isothermalisation},
since the heat distribution converges to a non-trivial isothermal state in all the space.
However, for dimension $N\geq 3$  uniqueness is lost in the
class of bounded solutions and some solutions decrease to zero as
$t\to\infty$, so that isothermalisation does not take place, see \cite{Kamin} and \cite{Eidus2}.

We are also facing the influence of the space dimension in the nonlocal case.
More precisely, if $J$ is compactly supported and $N=1,2$, we prove uniqueness
of bounded solutions. For dimension $N\geq3$, we need (as in the local case)
additional conditions on the behavior of $\rho$ at infinity to get uniqueness
of bounded solutions.

The case when $\rho$ is not integrable is also considered, which is more related to the study of the homogeneous case ($\rho\equiv 1$), see \cite{ChasseigneChavesRossi06} and \cite{BrandleChasseigneFerreira09}. For bounded solutions, the flux at infinity is so big that solutions go down to zero asymptotically while if the data is unbounded, the solution may go to infinity asymptotically as $t\to\infty$.

\

\noindent\textbf{Organisation and main results}

We first prove in Section \ref{sect:prelim} a comparison result which gives uniqueness for problem \eqref{eq:isot.0} in some class. In Section \ref{sect:ex.un} we study the existence in the class of bounded solutions,
which is obtained by approximation with Neuman problems in bounded domains. The main theorem is the following:
\begin{theorem}\label{thm:main.exist.uniq}
    Let $\rho>0$, continous, and $u_0$ be a bounded  nonnegative continuous function. Then
    there exits a nonnegative classical bounded solution of problem \eqref{eq:isot.0} such that
    $$\int_{\R^N}\rho(x)u(x,t)\d x=\int_{\R^N}\rho(x)u_0(x)\d x\,.$$
    Moreover, if either
    $$
    \rho(x)\ge \frac{\eta}{1+|x|^2}\,,
    $$
    or $N=1,2$ and $J$ has compact support then uniqueness holds in this class.
\end{theorem}

Section \ref{sect:u0.L1.rho} is
devoted to study the isothermalisation phenomenon for bounded initial data.
Then the main theorem is as follows:
\begin{theorem}\label{thm:main.isot}
    Let $\rho>0$, continuous, integrable and $u_0$ be a bounded  nonnegative continuous function. Then
    $u(x,t)\to\mathbb{E}_\rho(u_0)$ as $t\to
    \infty$ in $\L^p_{\rm loc}(\R^N)$ for any $1\leq p<\infty$; the convergence also
    holds in $\L^1(\rho)$:
    $$\lim_{t\to\infty}\int_{\R^N}|u(x,t)-\mathbb{E}_\rho(u_0)|\,\rho(x)\d x=0\,.$$
\end{theorem}
 If $\rho$ is not integrable but still $u_0\in\L^1(\rho)$
(that is, $\rho u_0\in\L^1(\mathbb R^N)$) the flux at infinity forces the solution to go to zero:
\begin{theorem}\label{thm:main.isot.0}
    Let $\rho>0$ and $u_0$ be a bounded nonnegative continuous function such that $u_0\in \L^1(\rho)$. If $\rho$ is not integrable in $\R^N$, then
    $u(x,t)\to 0$ as $t\to\infty$ in $\L^p_{\rm loc}(\R^N)$ for any $1\leq p<\infty$.
\end{theorem}
Finally, in Section \ref{sect.u0.not.L1.rho} we investigate the case of unbounded initial data.
We first prove an existence result in the class of unbounded solutions provided $\rho$ does not degenerate too rapidly at
infinity. More precisely, if
\begin{equation}\label{eq.decay}
    \rho(x)\geq\frac{\eta}{1+|x|^\gamma}\,\qquad \gamma\le 2\,,
\end{equation}
then we have an existence result for quadratic initial data:
\begin{theorem}\label{thm:main.nonisot}
    Let $u_0$ be a positive continuous function with at most quadratic growth
    at infinity. If the function $\rho$ satisfies \eqref{eq.decay}
    then there exists a solution of problem \eqref{eq:isot.0}.
\end{theorem}
More generally, if $\rho$ is integrable, we prove similar isothermalisation results for the minimal solution:
\begin{theorem}\label{thm:behaviour}
	We assume that $\rho$ is a continuous positive integrable function in $\R^N$ and that
  	$u_0$ is a continuous nonnegative function, possibly unbounded, such that there exists a solution $u$.
	Noting $\underline{u}$ the minimal solution, the following holds:

	i) If $u_0\in \L^1(\rho)$ the isothermalisation takes place in $\L^1(\rho)$,
	$$\lim_{t\to\infty}\int_{\R^N}|\underline{u}(x,t)-\mathbb{E}_\rho(u_0)|\,\rho(x)\d x=0\,.$$

	ii) If $u_0\notin \L^1(\rho)$, we have that for all $1\le p<\infty$,
	$$
	\lim_{t\to\infty} \underline{u}(x,t)= \infty \quad \mbox{in}\quad \L^p_{loc}(\mathbb R^N).
	$$
\end{theorem}

In the case of nonintegrable $\rho$'s with $u_0\not\in \L^1(\rho)$ the asymptotic behavior is more difficult to treat.
For instance, if $\rho\equiv 1$, the solutions
$$
u(x,t)=|x|^2+\mathbb{V}(J) t,\qquad \mbox{and}\qquad u(x,t)=1
$$
have different behavior.
Thus, there is a balance between $\rho$, $J$, and the initial data $u_0$ which is not easy to handle and the question remains open.

\

% ------------------------------------------------------------------------

\subsection*{Acknowledgment}
Both authors partially supported by project MTM2008-06326-C02-02 (Spain).
R.~Ferreira is also partially supported by grant GR58/08-Grupo 920894.
We would like to thank Jorge~Garcia-Meli\'an for giving us a proof that
$J$-harmonic bounded functions are constant (which was reproduced in Lemma~\ref{lem.strong.conv}). %------------------------------------------------------------------------

%%%%%%%%%%%%%%%%%%%%%%%%%%%%%%%%%%%%%%%%%%%%%%%%%%%%%%%%%%%%%%%%%%%%
%%%%%%%%%%%%%%%%%%%%%%%%%%%%%%%%%%%%%%%%%%%%%%%%%%%%%%%%%%%%%%%%%%%%
\section{Preliminaries}\label{sect:prelim}

Let us specify first what is the notion of solution that we use:
\begin{definition}\label{def:isot}
    Let $u_0\in \L^1_{\rm loc}(\R^N)$.
    By a stong solution of \eqref{eq:isot.0} we mean a function $u\in \C^0\big([0,\infty);\L^1_{\rm loc}(\R^N)\big)$
    such that $u_t,J\ast u\in \L^1_{\rm loc}(\R^N\times(0,\infty)\big)$, the equation
    is satisfied in the $\L^1_{\rm loc}$-sense and such that $u(x,0)=u_0(x)$
    almost everywhere in $\R^N$.
\end{definition}
We shall consider also solutions with more regularity:
\begin{definition}\label{def:isot.classical}
    A classical solution of \eqref{eq:isot.0} is a solution such that moreover
    $u,u_t,J\ast u\in \C^0(\R^N\times[0,\infty))$ and the equation holds in the classical sense
    everywhere in $\R^N\times[0,\infty)$.
\end{definition}
A classical sub or supersolution is defined as usual with inequalities instead of equalities
in the equation.

Now let us state a simple regularity Lemma, which contains a technical trick
that we shall use several times in the sequel:
\begin{lemma}\label{lem:reg}
    Let $u$ be a strong solution of \eqref{eq:isot.0}. We assume moreover that $u_0$ is continuous
    in $\R^N$ and that the convolution term $J\ast u$ is continous in $\R^N\times[0,\infty)$.
    Then $u$ and $u_t$ are also continuous in $\R^N\times[0,\infty)$
    and $u$ is a classical solution.
\end{lemma}
\begin{proof}
    We introduce the following transform:
    \begin{equation}\label{eq:transform}
        \mathcal{T}_\rho[u](x,t):=\e^{t/\rho(x)}u(x,t)\,.
    \end{equation}
    A straightforward calculus show that $v=\mathcal{T}_\rho[u]$ satisfies
    $$v_t=\frac{\e^{t/\rho(x)}}{\rho(x)}(J\ast u)(x,t)\,,$$
    which is a continuous function in $\R^N\times[0,\infty)$.
    Integrating between $0$ and $t$ we get:
    $$v(x,t)=\int_0^t \partial_t v(x,s)\d s+v(x,0)=\int_0^t \partial_t v(x,s)\d s+u(x,0)\,,$$
    hence $v$ is continuous in $\R^N\times[0,\infty)$. This implies that $u$ is also continuous
    in $\R^N\times[0,\infty)$, and the equation holds in the classical sense.
\end{proof}
\begin{remark}\label{rem:reg}
    {\rm It is well-known in the convolution theory that under one of the following assumptions,
    the convolution term is continuous:\\
    (i) $u$ is bounded (since $J$ is integrable);\\
    (ii) $J$ compactly supported and $u$ locally integrable.}
\end{remark}

The following lemma concerns the comparison of classical sub/supersolutions of the problem. In order to prove that we need to find a strict supersolution $\psi$, which satisfies
\begin{equation}\label{eq:def.psi}
    \begin{aligned}
    &\psi\in\mathrm{C}^0(\R^N\times[0,\infty))\,,\psi\geq0\,,\\\
    & \psi(x,t)\to+\infty\text{ as }|x|\to\infty\text{ uniformly for }t\in[0,\infty)\,,\\
    &\rho\psi_t> J\ast\psi-\psi\,.\\
    %&\psi(x,0)\ge u_0(x)
    \end{aligned}
\end{equation}

\begin{lemma}\label{lem.psi}In each of the following three cases, if $A>0$ and $\lambda$ is large enough, the function $\psi$ satisfies \eqref{eq:def.psi}:
\begin{itemize}
\item[(i)] If function $\rho$ satisfies
$$
\rho(x)\ge \frac{\eta}{1+|x|^2},
$$
then
$$
\psi(x)=A e^{\lambda t} (1+|x|^2)\,.
$$

\item[(ii)] If  $J$ has compact support and $N=1$,
$$
\psi(x)=A e^{\lambda t} (1+|x|)\,.
$$

\item[(iii)] If  $J$ has compact support and $N=2$,
$$
\psi(x)=A e^{\lambda t} \big(1+(\ln|x|)_+\big)\,.
$$
\end{itemize}
\end{lemma}

\begin{proof} The regularity and the behaviour at infinity of $\psi$ is clear, so we only need to prove that $\psi$ is a strictly supersolution.

In the first case a direct computation gives us that,
$$
J*\psi-\psi= A e^{\lambda t} \mathbb{V}(J)\,.
$$
On the other hand,
$$
\rho \psi_t= \rho A \lambda e^{\lambda t} (1+|x|^2)\ge A \lambda \eta e^{\lambda t}\,.
$$
Then, for $\lambda > \mathbb{V}(J)/\eta$, we are done.

In the other two cases, we notice first that $\psi$ is harmonic outside the origin.
Since $J$ is a probability density supported in $B(0,R)$, this implies that
$$
J*\psi-\psi=0,\qquad \mbox{for }\qquad |x|\ge 2R\,.
$$
On the other hand, since $\psi$ is continuous in $B(0,2R)$, there exists $K>0$ such that
$$
J*\psi-\psi \le K e^{\lambda t} \qquad \mbox{for }\qquad |x|<2R\,.
$$
Moreover, if $|x|\le 2R$ the function $\rho$ is bounded from below, so there exists a constant $\tilde K$ such that
$$
\rho \psi_t=\rho \lambda \psi(x,t)\ge \tilde K \lambda  e^{\lambda t}\,.
$$
Therefore, taking $\lambda\ge K/\tilde K$, we get that $\psi$ is a strictly supersolution.
\end{proof}

\begin{remark}In dimension $N\ge3$ and $J$ compactly supported, the function
$\psi(x,t)=e^{\lambda t}\min(1,|x|^{2-N})$ would give a strict supersolution,
however it goes to zero as $|x|\to\infty$. Hence it cannot be used to obtain
a comparison principle for bounded solutions.
\end{remark}

\begin{lemma}\label{lem:comp}
Let $\psi$ satisfy~\eqref{eq:def.psi}. Let  $\underline{u}$ be a classical subsolution of~\eqref{eq:isot.0} and $\bar u$ a
	classical supersolution of~\eqref{eq:isot.0} such that $\underline u(x,0)\leq \bar u(x,0)$. If for all  $t\ge 0$,
\begin{equation}\label{eq.condicion}
\limsup_{|x|\to\infty}\frac{\underline u-\overline u}{\psi}\ \le 0
\end{equation}
then $\overline{u}\geq \underline{u}$ in $\R^N\times\R_+$.
\end{lemma}
\begin{proof}
We consider the  function
$$
w_\delta=\underline u-\overline u- \delta\psi,
$$
which satisfies the inequality
\begin{equation}\label{eq:phi.delta}
	\rho(x)(w_\delta)_t<(J*w_\delta)(x,t)-w_\delta(x,t)\,.
\end{equation}
%The function $w_\delta$ is continuous in $\R^N\times[0,\infty)$ and notice that $w_\delta\geq\delta>0$ at time $t=0$.
%There are two options: \\
%\noindent $(i)$ either $w_\delta\ge \delta/2$ in $\R^N\times(0,\infty)$;\\
%\noindent $(ii)$ or there exists a point $(x,t)\in\R^N\times(0,\infty)$
%such that $w_\delta(x,t)< \delta/2$.
%
%Let us assume that we are in this second case and let $t_0$ be defined as follows:
%$$t_0:=\inf\{t>0:\exists x\in\R^N\,,\ w_\delta(x,t)<\delta/2\}<\infty\,.$$
%
%From \eqref{eq.condicion} we get that for $t\ge 0$
%\begin{equation}\label{eq.psi}
%\lim_{|x|\to\infty} w_\delta\ge \delta\,.
%\end{equation}
%Then, the only possible points where $w_\delta$ may reach a level above $\delta/2$ are located inside the fixed ball $B_{R_\delta}$.
%Thus, up to a finite number of terms, any minimizing sequence remains inside the compact set
%$\overline{B}_{R_\delta}\times[0,t_0+1]$. After extraction and using the continuity of $w_\delta$, we get that the inf is attained
%and there is a point $(x_0,t_0)$ such that $w_\delta(x_0,t_0)=\delta/2$. Of course this implies that $t_0$ cannot be zero, since $w_\delta\geq \delta$ at $t=0$.
%
%
%So, $t_0>0$ and necessarily $w_\delta\ge \delta/2$ in $\R^N\times(0,t_0)$,
%which implies also $J*w_\delta\ge \delta/2$ for $t\in(0,t_0)$, everywhere in $\mathbb{R}^N$. On the other hand, $(w_\delta)_t(x_0,t_0)\le 0$. Then, using \eqref{eq:phi.delta} we obtain the following contradiction

From \eqref{eq.condicion}, we deduce that on $\R^N\times[0,T]$, $w_\delta$ attains its maximum at a point $(x_0,t_0)$. Observe that at this point,
$$	\rho(x_0)(w_\delta)_t(x_0,t_0)<(J*w_\delta)(x_0,t_0)-w_\delta(x_0,t_0)\leq0$$
but if $t_0\in(0,T)$, $\partial_t w_\delta(x_0,t_0)=0$ and if $t_0=T$,
$\partial_t w_\delta(x_0,t_0)\geq0$. Hence in each of these cases we reach a contradiction. We are left with the last possibility, $t_0=0$, which implies $\max w_\delta\leq0$. Passing to the limit as $\delta\to0$, since $T>0$ is arbitrary we obtain the comparison
    $$
        \underline{u}-\overline{u}\le 0 \quad \mbox{for all} \quad x\in\R^N\,,\ t\ge 0\,.
    $$
\end{proof}

As direct corollary, we have the following uniqueness result,
which is valid only in a suitable class of solutions, as it is the case for the local heat equation.

\begin{corollary}
\label{thm:uniqueness}
    Let  $\psi$ satisfy~\eqref{eq:def.psi}. Let $u_0$ a continuous function which grows strictly less than $\psi$.
    Then there exists at most one classical solution $u$ of~\eqref{eq:isot.0}
    such that $|u| $ grows strictly less than $\psi$.
\end{corollary}

%%%%%%%%%%%%%%%%%%%%%%%%%%%%%%%%%%%%%%%%%%%%%%%%%%%%%%%%%%%%%%%%%%%%%%
%%%%%%%%%%%%%%%%%%%%%%%%%%%%%%%%%%%%%%%%%%%%%%%%%%%%%%%%%%%%%%%%%%%%%%

\section{Existence and uniqueness of bounded solutions}\label{sect:ex.un}

We obtain a existence result by approximation. We consider the following Neumann problem, where $\chi_n$ denotes the indicator of the ball $B_n$:

\begin{equation}\label{eq:neuman}
    \left\{\begin{array}{ll}
    \displaystyle
        \rho(x)\partial_t u_n = \int_{B_n}\big(u_n(y,t)-u_n(x,t)\big)J(x-y)\d y\,, \qquad & x\in B_n\,,\\
        u_n(x,0) = u_0\chi_n(x)\,,
    \end{array}\right.
\end{equation}

First, we observe that  this problem is possed in a bounded domain, so we have controlled the point where the maximum is attained. Therefore,
we can take $\psi=0$ in lemma \ref{lem:comp} to prove that a comparison principle holds.

%\begin{lemma}
%Let $\bar{u} \in \C^1([0,\infty);\C^0(\overline{B_n}))$ and $\underline{u} \in \C^1([0,\infty);\C^0(\overline{B_n}))$  be a supersolution and a subsolution of (\ref{eq:neuman}). Then $\bar u\ge \underbar u$.
%\end{lemma}
%\begin{proof}
%We consider the function
%$$
%w_{\delta} =\bar u-\underbar u +\delta e^t
%$$
%Which satisfies,
%\begin{equation}\label{eq.w}
%\rho(x) (w_{\delta})_t =\int_{B_n} J(x-y) (w_{\delta}(y,t)-w_{\delta}(x,t))\d y + \delta e^t \rho(x).
%\end{equation}
%The function $w_\delta$ is continuous in $\overline{B_n}\times[0,\infty)$ and notice that $w_\delta\ge \delta>0$ at time $t=0$. There are two options: \\
%\noindent $(i)$ either $w_\delta\geq \delta/2$ in $\overline{B_n}\times(0,\infty)$;\\
%\noindent $(ii)$ or there exists a point $(x,t)\in\overline{B_n}\times(0,\infty)$
%such that $w_\delta(x,t)<\delta/2$.
%
%Let us assume that we are in this second case and let $t_0$ be defined as follows:
%$$t_0:=\inf\{t>0:\exists x\in\overline{B_n}\,,\ w_\delta(x,t)<\delta/2\}<\infty\,.$$
%
%
%
%Since the function $w_\delta$ is continuous, we get that the inf is attained
%and there is a point $(x_0,t_0)$ such that $w_\delta(x_0,t_0)=\delta/2$. Of course this implies that $t_0$ cannot be zero, since $w_\delta\geq \delta$ at $t=0$.
%
%
%So, $t_0>0$ and necessarily $w_\delta\geq \delta/2$ in $\overline{B_n} \times(0,t_0)$ and $(w_\delta)_t(x_0,t_0)\le 0$. Then, using \eqref{eq.w} we obtain the following contradiction
%$$
%0\ge \rho (w_\delta)_t(x_0,t_0)= \int_{B_n} J(x_0-y) (w_{\delta}(y,t_0)-w_{\delta}(x_0,t_0))\d y + \delta e^{t_0} \rho(x_0)>0\,.
%$$
%\end{proof}

\begin{lemma}\label{lem:neuman}
    Let $\rho>0$ be a continuous function in $\R^N$ and $u_0$ be a nonnegative bounded continuous function.
    Then for any $n\geq0$, there exists a unique solution $u_n\in\C^1([0,\infty);\C^0(\overline{B_n}))$, which satisfies $0\le u_n\le \|u_0\|_\infty$. Moreover, the following conservation law holds:
    $$\int_{B_n} u_n(x,t)\rho(x)\d x = \int_{B_n} u_0(x)\rho(x)\d x\,.$$
\end{lemma}
\begin{proof}
    Following \cite{ChasseigneChavesRossi06}, we consider $t_0>0$ to be fixed later on and the Banach space
    $X_{t_0}=\C^0([0,t_0]\times \overline{B_n})$ equipped with the
    norm $\vertt w \vertt = \max\big\{\Vert w(\cdot,t)\Vert_{\L^\infty(\overline{B_n})}\,, 0\leq t\leq t_0\big\}$. Then we define an operator
    ${T}:X_{t_0}\to X_{t_0}$ as follows:
    $${T}_{w_0}(w)(x,t):=w_0(x)+\frac{1}{\rho}\int_0^{t}\int_{B_n}J(x-y) (w(y,s) - w(x,s))\d y\d s$$
    A straightforward calculus, using that $\rho\geq \rho_0$ in $B_n$ shows that:
    $$
    \vertt{T}_{w_0}(w)-{T}_{z_0}(z)\vertt\leq\| w_0-z_0\|_{\L^\infty(\overline{B_n})}+\frac{2 t_0}{\rho_0}\vertt w-z\vertt\,.
    $$
    We deduce that for any fixed initial data $u_0\in \L^\infty(\overline{B_n})$, if $t_0>0$ is sufficiently small (say $t_0<\rho_0/4$),
    ${T}_{u_0}$ is a strict contraction in the Banach space $X_{t_0}$. Hence there exists a unique function $u_n\in X_{t_0}$
    such that $T_{u_0}(u_n)=u_n$, which means that for any $x\in \overline{B_n}$ and $t\in[0,t_0]$, one has
    $$
    u_n(x,t)=u_0(x)+\frac{1}{\rho}\int_0^{t}\int_{B_n}J(x-y) (u_n(y,s) - u_n(x,s))\d y\d s\,.
    $$
    This in particular implies that $u(0)=u_0$ and moreover, that the weak derivative $\partial_t u$ exists, which is given by
    $$\partial_t u_n = \frac{1}{\rho}\int_{B_n}J(x-y) (u_n(y,s) - u_n(x,s))\d y \,.$$
    Since the right hand side of the equation is a continuous function in $\overline{B_n}\times[0,t_0]$,
    we have that $\partial_t u$ is also continuous in $\overline{B_n}\times[0,t_0]$.  to extend the solution to $[0,\infty)$ we take as initial data $u(x,t_0)\in C^0(\overline{B_n})$ and obtain a solution up $[0,2t_0]$. We then iterate the procedure to construct a solution for all time $t>0$.

    On the other hand, we note that the constants functions $\bar u=\|u_0\|_\infty$ and $\underline u=0$ are a supersolution and subsolution respectively, then by comparison
    $$
    0\le u_n\le \|u_0\|_\infty\,.
    $$

    Finally, we observe that
    $$
    \rho(x) u_n(x,t)-\rho(x) u_0(x)=\int_0^{t}\int_{B_n}J(x-y) (u_n(y,s) - u_n(x,s))\d y\d s\,.
    $$
    Then, integrating in $x$ and apply Fubini's theorem in the right hand side of the equation, we get de desired conservation law.
\end{proof}

\begin{proposition}\label{prop:limit.neuman}
    Let $\rho>0$ be a continuous function in $\R^N$
    and $u_0$ be a nonnegative continuous bounded function and
    for any integer $n>0$, let $u_n$ be the solution constructed in Lemma \ref{lem:neuman}. Then
    as $n\to\infty$, along a subsequence $u_n\to u$ in $\L^1_{\rm loc}(\R^N\times[0,\infty))$
    where $u$ is a nonnegative classical bounded solution of \eqref{eq:isot.0}.
\end{proposition}

\begin{proof}
    Since the sequence $u_n$ is bounded, it converges (along a subsequence still denoted $u_n$)
    in $\L^\infty$-weak$\ast$ to some nonnegative and
    bounded function $u$. Since $J\in\L^1$, this implies that $J\ast u_n\to J\ast u$ strongly
    and that $\partial_t u_n\to \partial_t u$ in the sense of distributions. We can then pass to the
    limit in the sense of distributions but here also we want a better convergence.

    Let us introduce a modified version of transform $\mathcal{T}$ as follows:
    $$v_n(x,t):=\e^{t(J\ast\chi_n)(x)/\rho(x)}w(x,t)\,.$$
    Since \eqref{eq:neuman} can be written as
    $$\rho(x)\partial_t(u_n)(x,t)=\big[J\ast (u_n\chi_n)\big](x,t) - (J\ast\chi_n)(x) u_n(x)\,,$$
    it follows immediately that $v_n$ satisfies the equation
    $$\partial_t v_n = \frac{\e^{t(J\ast\chi_n)(x)/\rho(x)}}{\rho(x)}J\ast u_n\,.$$

    This implies that $\partial_t v_n$ converges strongly on compact sets of $\R^N\times[0,\infty)$,
    and so does $v_n(x,t)=u_0(x)\chi_n(x)+\int_0^t (v_n)_t(x,s)\d s$. Then $u_n$ also converges strongly
    on compact sets of $\R^N\times[0,\infty)$ to its limit $u$.

    Passing to the limit in the equation, we see that $u$ is a strong solution of \eqref{eq:isot.0},
    which implies that it is a classical solution of this equation -- this follows from Lemma
    \ref{lem:reg}.
    %
%    Any other converging subsequence
%    leads to the same solution $u$ so that all the sequence $u_n$ converges to $u$ in
%    $\L^1_{\rm loc}(\R^N\times[0,\infty))$.
\end{proof}

\begin{proof}[Proof of Theorem \ref{thm:main.exist.uniq}]
It is just the combination of Proposition \ref{prop:limit.neuman}, Corollary \ref{thm:uniqueness} and Remark \ref{lem.psi}.
\end{proof}

Now, we pass to the limit in the conservation low to obtain:

\begin{theorem}\label{prop.cons.law.rho}
    Let $\rho>0$, continuous and $u_0\in\C^0(\R^N)\cap\L^\infty(\R^N)\cap\L^1(\rho)$.
    Then \\
    $(i)$ for any $t>0$, $u(\cdot,t)\in\L^1(\rho)$ and
    $$
    \int_{\mathbb R^N} \rho(x) u(x,t)\d x\le \int_{\mathbb R^N} \rho(x) u_0(x)\d x\,.
    $$
    $(ii)$ if moreover we assume that $\rho\in\L^1(\R^N)$ then the conservation law holds:
    $$
    \int_{\mathbb R^N} \rho(x) u(x,t)\d x = \int_{\mathbb R^N} \rho(x) u_0(x)\d x\,.
    $$
\end{theorem}
\begin{proof}
    Since $u_0\in\L^\infty\cap\mathrm{C}^0$, the approximating sequence $\{u_n\}$ of Proposition~\ref{prop:limit.neuman}
    satisfies that (up to a subsequence) $u_n\to u$ almost everywhere, with $0\leq u_n\leq \|u_0\|$ and
    \begin{equation}\label{eq.conser_n}
    \int_{\R^N} u_n(x,t)\chi_n(x)\rho(x)\d x = \int_{\R^N} u_0(x)\chi_n(x)\rho(x)\d
    x\,.
    \end{equation}
Since  $u_0\in\L^1(\rho)$, the dominated convergence theorem yields
the convergence of the right-hand side integral as $n\to\infty$.
For the left hand side, using Fatou's Lemma we obtain that
    $$
    \int_{\mathbb R^N} \rho(x) u(x,t)\d x\le \liminf_{n\to\infty}
    \int_{\mathbb R^N} \rho(x) u_n(x,t)\chi_n(x)\d x =
    \int_{\mathbb R^N} \rho(x) u_0(x)\d x\,,
    $$
    which proves assumption $(i)$.

    Finally, if we assume that $\rho\in\L^1(\R^N)$ we can use also the dominated convergence theorem for the sequence
    $u_n\chi_n\rho$, which is bounded by $\|u_0\|_\infty\,\rho\in\L^1(\R^N)$. We then pass to the limit in
    the left-hand side of \eqref{eq.conser_n} and get $(ii)$.
\end{proof}

%\begin{proposition}
%   For any $u_0\in\L^1(\rho)\cap\L^\infty(\R^N)$, there exists a nondecreasing function $\ell(t)\geq0$ such that
%   $\ell(0)=0$ and $\lim \ell(t)=\ell_\infty <\infty$ and
%   $$
%       \int_{\R^N}\rho(x)u(x,t)\d x + \ell(t) = \int_{\R^N}\rho(x)u(x,t)\d x\,.
%   $$
%\end{proposition}
%
%\begin{proof}
%   First, it is clear that for any $t\geq0$, such a $\ell(t)$ exists, since it can be defined by:
%   $$\ell(t):=\int \rho u_0-\int \rho u(t)\,.$$
%   By definition, $0=\ell(0)\leq\ell(\cdot)\leq\int \rho u_0$. Now we claim that it is noondecreasing:
%   indeed, since $t\mapsto u(t+s)$ may be viewed as a solution with initial data $u(s)$, then it follows that
%   for any $s,t>0$,
%   $$\int \rho u(t+s)\leq \int \rho u(s)\Rightarrow \ell(t+s)\geq\ell(s)\,.$$
%   So, the limit $\ell_\infty=\lim \ell(t)$ as $t\to\infty$ is also well-defined.
%\end{proof}
%
%\begin{remark}{\rm
%$(i)$ this will be enough to conclude that the $\omega$-limit set is a singleton, since we will get, if $\rho$ is integrable:
%$$\omega(u_0)=\Big\{ \mathbb{E}_\rho(u_0)-\frac{\ell_\infty}{\int \rho} \Big\}$$
%$(ii)$ it seems difficult for the moment to find conditions that insure $\ell(t)=0$;
%at least we know that if $\rho$ is not integrable, the flux at infinity is so big that
%$\ell_\infty=\int \rho u_0$ and finally, $\int\rho u(t)\to0$.}
%\end{remark}

%%%%%%%%%%%%%%%%%%%%%%%%%%%%%%%%%%%%%%%%%%%%%%%%%%%%%%%%%%%%%%%%%%%%
%%%%%%%%%%%%%%%%%%%%%%%%%%%%%%%%%%%%%%%%%%%%%%%%%%%%%%%%%%%%%%%%%%%%

\section{Asymptotic behaviour for bounded solutions}\label{sect:u0.L1.rho}

We shall now derive our main results concerning the asymptotic behaviour for \eqref{eq:isot.0}.
We divide the proof in several steps.

%%%%%%%%%%%%%%%%%%%%%%%%%%%%%
\subsection{Weak limit}
\label{subsec.wl}

This first step is easy, it only comes from the fact that the solution is globally bounded:
\begin{lemma}\label{lem.weak.conv}
    Let $\rho>0$, continuous and $u_0\in\C^0(\R^N)\cap\L^\infty(\R^N)$. Let $u$ be
    a bounded classical solution.
    Then for any $s>0$ there exists a subsequence $t_k\to+\infty$ such that
    the following limit exists in $\L^\infty$-weak*:
    $$u_\infty(x,s):=\lim_{t_k\to\infty}u(x,s+t_k)\,.$$
\end{lemma}
\begin{proof}
    Since $u$ is bounded, thus there exists a
    subsequence $t_k\to\infty$ such that $u(\cdot,s+t_k)$ converges in $\L^\infty$-weak* to
    a function $u_\infty(\cdot,s)\in\L^\infty(\R^N)$.
\end{proof}
%%%%%%%%%%%%%%%%%%%%%%%%%%%%%
\subsection{Lyapounov functional}
\label{subsec.lf}

We now want a stronger result, so we use a Lyapounov functional:
\begin{lemma}\label{lem.cota.lyapunov}
Let $t_0>0$. Assume the hypotheses of Lemma \ref{lem.weak.conv} and that $u_0\in\L^2(\rho)$. Then there exists a constant $C=C(u_0,\rho,t_0)$ such that
for all $t\ge t_0>0$,
  $$
    \int_t^\infty  \int_{\mathbb R^N} \rho(x) (u_t)^2(x,s) \d s \le C.
  $$
\end{lemma}

\begin{proof}
    We first go back to the approximating scheme (\ref{eq:neuman}). Notice that the following quantity
    $$
    F[u_n](t)=\int_{B_n}\int_{B_n}
    J(x-y)(u_n(x,t)-u_n(y,t))^2\d x\d y
    $$
    is a Lyapunov functional which is nonincreasing along the evolution orbits at the level $n$. Indeed, multiplying the equation by $u$ and integrating in $B_n$ we get
    $$\begin{aligned}
        \int\rho(x)\partial_t u_n\cdot u_n\d x &= \iint J(x-y)u_n(x)u_n(y)\d x\d y - \int u_n(x)^2\d x \\
        &= \frac{1}{2}\iint J(x-y)\big(2u_n(x)u_n(y)-2u_n(x)^2\big)\d x\d y\\
        &= \frac{1}{2}\iint J(x-y)\big(2u_n(x)u_n(y)-u_n(x)^2-u_n(y)^2\big)\d x\d y\\
        &= -\frac{1}{2}F[u_n](t)\,.\end{aligned}$$
    Therefore,
    \begin{equation}\label{eq.lyapounov.int}
    F[u_n](t)=- \frac{\d}{\d t} \int_{B_n}  \rho_n(x) u_n^2(x,t)\d x\,.
    \end{equation}
    In a similar way, multiplying the equation by $\partial_t u_n$ and integrating in space we obtain
    \begin{equation}\label{eq.lyapounov.der}
    \frac{\d}{\d t}F[u_n](t)=- 4 \int_{B_n} \rho(x)((u_n)_t)^2(x,t)\d x\,.
    \end{equation}

    Integrating  \eqref{eq.lyapounov.int} we have that for some $C'=C'(u_0,\rho)$,
    $$
    \int_0^t F[u_n](s)\d s = 2\int_{B_n}  \rho(x) u_0^2(x) \chi_n(x) \d x-2
    \int_{B_n}  \rho(x) u_n^2(x,t)\d x\le C'\,.
    $$
    Indeed, as $u_0\in\L^2(\rho)$ we have by monotone convergence that
    $$
    \int_{\mathbb{R}^N}  \rho(x) u_0^2(x) \chi_n(x) \d x \to \int_{\mathbb{R}^N}  \rho(x) u_0^2(x) \d x<\infty\,.
    $$
    Using now the monotonicity in $t$ of $F[u_n](t)$, we get
    $$
    tF[u_n](t)\le \int_0^t F[u_n](s)\d s\le C'\,.
    $$
    Moreover, $F[u_n](t)$ is
    positive so that, integrating \eqref{eq.lyapounov.der}, we get for any $t\ge t_0$:
    $$
    \int_t^\infty \int_{\mathbb{R}^N}\rho ((u_n)_t)^2(x,s)\d x\d s\leqslant \frac14\,F[u_n](t)\le \frac{C'}{4t_0}.
    $$
    Using Fatou's Lemma and the fact that $\rho (u_n)_t$ converges strongly to $\rho u_t$, we obtain
    the desired result.
\end{proof}

As and immediate consequence of this result we obtain:

\begin{lemma}\label{lem.lyapounov}
  Assume the hypotheses of Lemma \ref{lem.weak.conv} and that $u_0\in\L^2(\rho)$. For all sequence $t_k\to \infty$ and $s>0$,
$$
\|\sqrt{\rho(\cdot)} \,u(\cdot,s+t_k)-\sqrt{\rho(\cdot)}\, u(\cdot,t_k)\|_{\L^2(\mathbb{R}^N)}^2
\to 0 \quad \mbox{as }n\to\infty\,.
$$
Hence, the limit function $u_\infty(x,s)$ does not depend on the variable $s>0$.
\end{lemma}
\begin{proof}
    Notice that for all sequence $t_k\to\infty$, we get
    $$
    \begin{aligned}
    \big\|\sqrt{\rho(\cdot)} \,u(\cdot,s+t_k)-\sqrt{\rho(\cdot)}\, u(\cdot,t_k)\big\|_{\L^2(\mathbb{R}^N)}^2 &
    = {\ds \int_{\mathbb{R}^N}\rho(x)\Big(\int_{t_k}^{t_k+s} u_t(x,\sigma)\d \sigma}\Big)^2\d x\\
    & \le s{\ds \int_{\mathbb{R}^N}\int_{t_k}^{t_k+s}\rho(x)(u_t)^2(x,\sigma)\d \sigma \d x},
    \end{aligned}
    $$
    which goes to zero as $t_k\to+\infty$.
\end{proof}

%%%%%%%%%%%%%%%%%%%%%%%%%%%%%

\subsection{The $\omega-$limit set}
\label{subsec.sl}

We define the $\omega-$limit set as follows
$$
\omega(u_0)=\{u_\infty\in \mathrm{C}^0(\R^N)\,:\,\exists t_j\to\infty \mbox{ such that } u(\cdot,t_j)\to u_\infty(\cdot)
\mbox{ in $\L^\infty$-weak*}  \}
$$

\begin{lemma}\label{lem.strong.conv}
    Under the hypotheses of Lemma \ref{lem.weak.conv} and that $u_0\in\L^2(\rho)$, the $\omega$-limit set is reduced to constants.
\end{lemma}
\begin{proof}
    Since $u(x,s+t_k)$ converges weakly in $\L^\infty$-weak*, then as $t_k\to\infty$,
    $$(J\ast u)(x,s+t_k)=\int J(x-y)u(y,s+t_k)\d y\stackrel{\rm pointwise}{\longrightarrow} (J\ast u_\infty)(x,s)\,.$$
    Moreover, since $u$ is bounded, the convergence of $J\ast u$ is also strong in $\L^1_{\rm loc}$.
    On the other hand, $\partial_t\rho u(s+t_k)\to \rho \partial_s u_\infty(s)$ in the sense of
    distributions. We then pass to the limit in the sense of distributions in the equation and get
    $$
    \rho(x)\frac{\partial}{\partial s} u_\infty(x,s) = J\ast u_\infty(x,s) -u_\infty(x,s)\,.
    $$

    Using Lemma \ref{lem.lyapounov} we know that $u_\infty$ is independent of $s$ so that
    $u_\infty$ is a bounded solution (in the sense of distributions) of
    $$J\ast u_\infty-u_\infty=0\ \text{ in }\ \R^N\,.$$
    We first deduce that $u_\infty$ is continuous because the convolution term is continous.
    Then we consider the Fourier transform in the space of tempered distributions:
    $$\hat{J}\cdot\hat{u}_\infty-\hat{u}_\infty=0\,.$$
    Since $|\hat{J}|<1$ except for $\hat{J}(0)=1$, we see that
    the (generalized) Fourier transform of $u_\infty$ has a support contained in $\{0\}$.
    This implies that $u_\infty$ is a polynom but since it is bounded, it has to be a constant.
\end{proof}

%%%%%%%%%%%%%%%%%%%%%%%%%
\subsection{Identification of the limit}
\label{subsec.cl}

We are now ready to identify the $\omega$-limit set.
\begin{lemma}\label{lem:identify}
We assume the hypotheses of Lemma \ref{lem.weak.conv} and that $u_0\in\L^1(\rho)$.
Then the following holds:\\[2mm]
$(i)$ if $\rho\in\L^1(\R^N)$, $\omega(u_0)=\big\{\mathbb{E}_\rho(u_0)$\big\}\,;\\[2mm]
$(ii)$ if $\rho\notin\L^1(\R^N)$, $\omega(u_0)=\big\{0\big\}$\,.
\end{lemma}
\begin{proof}
	Notice first that since $u_0\in\L^1(\rho)\cap\L^\infty$, then $u_0\in\L^2(\rho)$,
	hence we may use Lemma \ref{lem.strong.conv}. In the integrable case, $\rho\in\L^1$,
	we observe that as $u$ is uniformly bounded, the dominated convergence Theorem gives
    $$
    \int_{\mathbb R^N} \rho(x) u_n(x,s+t_j)\d x \to u_\infty \int_{\mathbb R^N} \rho(x) \d x\,.
    $$
    Therefore, by Theorem \ref{prop.cons.law.rho}-$(ii)$ we obtain
    $u_\infty=\mathbb{E}_\rho(u_0)$, so that the $\omega$-limit set is reduced to
    $\big\{\mathbb{E}_\rho(u_0)\big\}$.

    In the case $\rho\not\in\L^1(\R^N)$, we take a compact set $K$
    such that
    $$
    \int_K \rho(x)\d x>\int_{\R^N} \rho(x)u_0(x)\d x\,,
    $$
    which is always possible since $u_0\in\L^1(\rho)$.
    Using Theorem~\ref{prop.cons.law.rho}-$(i)$,
    Lemma~\ref{lem.strong.conv} and Fatou's Lemma, we obtain
    $$
    \begin{array}{rl}
    u_\infty \int_{K} \rho(x) \d x\le &\liminf\limits_{n\to\infty}
    \int_{K} \rho(x) u_n(x,s+t_j)\d x\\ \le &\liminf\limits_{n\to\infty}
    \int_{\mathbb R^N} \rho(x) u_n(x,s+t_j)\d x \\
    \le &\int_{\mathbb R^N} \rho(x) u_0(x)\d x\,.
    \end{array}
    $$
    This implies that necessarily $u_\infty=0$, hence the $\omega$-limit set is reduced to
    $\big\{0\big\}$.
\end{proof}

%%%%%%%%%%%%%%%%%%%%%%%%%%%%%
\subsection{Proofs of Theorems \ref{thm:main.isot} and \ref{thm:main.isot.0}}
\label{subsec.teos}

As consequence of the fact that the $\omega-$limit is given by only one function we can pass to the limit in the time variable without extracting any subsequence. Then it only remains to check that the convergence is better, which is done by using transform $\mathcal{T}_\rho$.

\begin{proof}[Proof of Theorem \ref{thm:main.isot.0}]
	Under the hypotheses of the Theorem, we have that $u_0$ is continuous, bounded,
	and $\rho$ not integrable, but nevertheless $u_0\in\L^1(\rho)\cap\L^2(\rho)$. Thus Lemmas  \ref{lem.weak.conv}
	and \ref{lem:identify} imply that for any $s>0$, at least along a subsequence $t_n\to\infty$ we have
    $u(x,s+t_n)\to0$ in $\L^\infty$-weak*. But the same arguments are valid for any other
    subsequence such that $u(x,s+{t'}_n)$ converges weakly. Since the
    limit is always zero, we deduce that for any $s>0$,
	 $$u(x,s+t)\mathop{\xrightarrow{\hspace*{1.5cm}}}_{t\to\infty}^{\L^\infty\text{-weak*}} 0\,,$$
	which implies that $\big(J\ast u(s+t)\big)$
	converges strongly in $\L^1_{\rm loc}$ as $t\to\infty$. Then,
	$$\rho(x)\partial_s u(x,s+t)=\big(J\ast u(s+t)\big)(x) - u(x,s+t)
	\mathop{\xrightarrow{\hspace*{1.5cm}}}_{t\to\infty}^{\L^\infty\text{-weak*}} 0\,.$$
	Even more, since $t\to+\infty$ we may assume that $t>t_0$ for some $t_0>0$
	and from lemma \ref{lem.cota.lyapunov} we obtain that for any compact set $K$,
$$\begin{aligned}
\int_t^\infty \left(\int_K \rho(x)|u_t|(x,s)\d x\right)^2\d s\le &
\int_t^\infty \left(\int_K \rho(x)|u_t|^2(x,s)\d x\right)\\
\le & \left(\int_K \rho(x) \d x\right)\d s \le C(K,u_0,\rho).\end{aligned}
$$
Then, at least for some sequence $t_k\to+\infty$, we have $\rho(x) \partial_s u(x,s+t_k)\to 0$ in $\L^1_{\rm loc}$. Summing up, we obtain that
    $$\lim_{t_k\to+\infty}u(x,s+t_k)=0\text{ in }\L^1_{\rm loc}\,.$$
	Of course, if $t\mapsto u(x,s+t)$ were to converge in $\L^1_{\rm loc}$ along another subsequence $t'_k\to\infty$, the limit
	would necessarily be zero, so that finally $u(\cdot,t)\to 0$ in $\L^1_{\rm loc}(\R^N)$ as $t\to\infty$.
    Moreover, since $t\mapsto u(\cdot,t)$ remains bounded in $\L^\infty(\R^N)$,
    we deduce that the convergence holds in $\L^p_{\rm loc}(\R^N)$ for any $1\le p<\infty$.
\end{proof}

\begin{proof}[Proof of Theorem \ref{thm:main.isot}]
	The first part is done exactly as in the proof of Theorem \ref{thm:main.isot.0}, except that $\rho$ is
	integrable here so that the limit is not zero, but $\mathbb{E}_\rho(u_0)$.
	To end the proof in this case, it only remains to prove
	the $\L^1(\rho)$ convergence. We fix $\varepsilon>0$ and choose $R>0$ big enough so that
    (remember that $\rho$ is integrable):
    $$\int_{|x|>R}\rho(x)\d x\leq\varepsilon\,.$$
    Then
    $$\int_{\R^N}|u(x,t)-u_\infty|\rho(x)\d x\leq 2\varepsilon\Vert u_0\Vert_\infty+\int_{|x|\leq R}
    |u(x,t)-u_\infty|\rho(x)\d x\,,$$
    and using the $\L^1_{\rm loc}$ convergence we get:
    $$\limsup_{t\to\infty}\int_{\R^N}|u(x,t)-u_\infty|\rho(x)\d x\leq 2\varepsilon\Vert u_0\Vert_\infty\,.$$
    Since $\varepsilon$ is arbitrary, we get that the limit is zero, which ends the proof.
\end{proof}

%%%%%%%%%%%%%%%%%%%%%%%%%%%%%%%%%%%%%%%%%%%%%%%%%%%%%%%%%%%%%%%%%%%%
%%%%%%%%%%%%%%%%%%%%%%%%%%%%%%%%%%%%%%%%%%%%%%%%%%%%%%%%%%%%%%%%%%%%
\section{Unbounded solutions}\label{sect.u0.not.L1.rho}

In this section we derive some results for unbounded initial data and solutions.
Let us mention that in the case $\rho\equiv1$, further results are to be found in
\cite{BrandleChasseigneFerreira09}. But here we still face the problem of the
space inhomogeneity implied by $\rho$.

%%%%%%%%%%%%%%%%%%%%%%%%%%%%%
\subsection{An existence result for unbounded solutions}

\begin{proposition}\label{prop:ex.unbounded}
    Let us assume that \eqref{eq:def.psi} holds.
    Then for any nonnegative $u_0\in \C^0(\R^N)$, satisfying $u_0(x)\leq \psi(x,0)$
    there exist a classical solution $u$ of \eqref{eq:isot.0} with $u(x,0)=u_0(x)$.
\end{proposition}

\begin{proof}
    Let us first consider an approximation $u_{0n}=u_0\cdot\chi_n$
    where $\chi_n$ is smooth, nonnegative, compactly supported and
    $\chi_{n}\nearrow1$. Let $u_n$ be the solution of \eqref{eq:isot.0} with initial
    data $u_{0n}$ given by Proposition \ref{prop:limit.neuman}, then by applying
    the comparison result for bounded solutions, the sequence $u_n$ is nondecreasing.

    Notice that $\psi$ is not bounded, but this is allowed in
    Lemma \ref{lem:comp}, which gives:
    $$ u_n(x,t)\leq \psi(x,t)\,.$$
    Hence the sequence $u_n$ converges to some $u$ and we are able to pass to the limit in $J\ast u_n$ by dominated convergence, using that $\psi$ is integrable with respect    to translations of $J$.

    Using now Lemma \ref{lem:reg}, we deduce that $u$ is a classical solution of \eqref{eq:isot.0} and the initial data of $u$ is $u_0$.
\end{proof}

\begin{remark}
    {\rm This construction does in fact give a minimal solution: if $u_1$ is any other solution, then it can
    be used as a supersolution for any $u_n$ and passing to the limit shows that $u\leq u_1$. One can think that if we restrict the initial data to grow strictly less than  $\psi$, then uniqueness holds because the comparison argument is valid in this class. However, it is not clear whether the constructed solution enters this class unless we know more about $u_0$, see~\cite{BrandleChasseigneFerreira09}.}
\end{remark}

%%%%%%%%%%%%%%%%%%%%%%%%%%%%%%%%%%%%%%%%%%%%%%%%%%%%%%%%%%%%%
\subsection{Asymptotic behaviour for unbounded solutions when $\rho$ is integrable}

We prove now that if $u_0$ is integrable with respect to $\rho$, the isothermalisation phenomenon occurs (whether infinite or not).
Notice that we gave sufficient conditions for existence of a minimal solution in the previous section. The first result is the following:

\begin{proposition}\label{prop:finite.isot}
    Let $u_0\in \C^0(\R^N)\cap \L^1(\rho)$, $\rho\in \L^1(\R^N)$ and assume there exists a solution $u$ such
    that $u(x,0)=u_0(x)$. Then, if $\underline{u}$ denotes the minimal solution, we have
    $$\lim_{t\to\infty}\int_{\R^N}|\underline{u}(x,t)-\mathbb{E}_\rho(u_0)|\,\rho(x)\d x=0\,.$$
\end{proposition}

\begin{proof}
Let $\underline{u}$ be the minimal solution and let us use the same monotone approximations that were used
in Proposition \ref{prop:ex.unbounded}. Since $u_n$ is bounded, Lemma \ref{lem:comp} implies that we have a bound
from above:
$$
u(x,t)\geq u_n(x,t)\text{ in }\R^N\times[0,\infty)\,.
$$
But since $u_n(x,0) \in \L^{\infty}(\mathbb R^N)$, Theorem \ref{prop.cons.law.rho} implies
$$
\int_{\mathbb R^N} u_n(x,t) \rho(x)\d x=\int_{\mathbb R^N} u_n(x,0) \rho(x)\d x.
$$
Moreover, the convergence of $u_n$ to $\underline{u}$ is monotone, so we can pass to the limit in the above equation
to obtain
$$
\int_{\mathbb R^N} \underline{u}(x,t) \rho(x)\d x=\int_{\mathbb R^N} \underline{u}(x,0) \rho(x)\d x.
$$
Using the above three equations we get that
$$
\int_{\mathbb R^N} |\underline{u}(x,t)- \mathbb{E}_\rho(u_0)|\rho(x)\d x \le I_1+I_2+I_3,
$$
where,
$$
I_1=\int_{\mathbb R^N} |\underline{u}(x,t)- u_n(x,t)|\rho(x)\d x =\int_{\mathbb R^N} (\underline{u}(x,0)- u_n(x,0))\rho(x)\d x\,.
$$
$$
I_2=\int_{\mathbb R^N} |u_n(x,t)- \mathbb{E}_\rho(u_0\chi_n)|\rho(x)\d x,
$$
$$
I_3=\int_{\mathbb R^N} |\mathbb{E}_\rho(u_0\chi_n)-\mathbb{E}_\rho(u_0)|\rho(x)\d x.
$$
Observe that $I_1$ and $I_3$ are independents of $t$ and tend to zero as $n\to \infty$. Moreover, $u_n$ satisfies the hypothesis of Theorem \ref{thm:main.isot} so that $I_2$ tends to zero as $t\to \infty$. Therefore, we first have
$$\limsup_{t\to+\infty}\int_{\mathbb R^N} |\underline{u}(x,t)- \mathbb{E}_\rho(u_0)|\rho(x)\d x \le I_1 + I_3\,,$$
so that taking the limit as $n\to\infty$ yields the desired result.
\end{proof}

In the case when $u_0\notin \L^1(\rho)$, then infinite isothermalisation
occurs:

\begin{proposition}\label{prop:infinite.isot}
    Let $\rho\in \L^1(\R^N)$ and $u_0\in \C^0(\R^N)$ such that $u_0\notin \L^1(\rho)$. Then for any solution $u$ with initial data $u_0$
	and any $1\leq p<\infty$, the following asymptotic behaviour holds:
    $$\lim_{t\to+\infty}u(x,t)=+\infty\text{ in } \L^p_{\rm loc}(\R^N)\,.$$
\end{proposition}

\begin{proof}
    As before, if there exists a solution, then we can approximate the minimal solution $\underline{u}$
    by the family $u_n$ used
    in Proposition \ref{prop:ex.unbounded}. Since this approximation is monotone,
    $$\underline{u}(x,t)\geq u_n(x,t)\text{ in }\R^N\times[0,\infty)\,.$$
    But $u_n(x,0)$ satisfies the hypotheses of Theorem \ref{thm:main.isot} so that
    $$\liminf_{t\to+\infty}\underline{u}(x,t)\geq\lim_{t\to+\infty}u_n(x,t)=c_n\,,$$
    where $c_n=\mathbb{E}_\rho(u_n(x,0))$, the limit holding in all $\L^p_{\rm loc}(\R^N)$.
    Hence passing to the limit as $n\to+\infty$, we obtain the result for $\underline{u}$ since
    $c_n\to\mathbb{E}_\rho(u_0)=+\infty$, thus the same holds for any other solution.
\end{proof}

Theorems \ref{thm:main.nonisot} and \ref{thm:behaviour} follow from the conjunction of Propositions and \ref{prop:ex.unbounded}, \ref{prop:finite.isot} and \ref{prop:infinite.isot}.

%%%%%%%%%%%%%%%%%%%%%%%%%%%%%%%%%%%%%%%%%%%%%%%%%%%%%%%%%%%%%%%%%%
%\subsection{Asymptotic behaviour for unbounded solutions when $\rho$ is not integrable}
%\label{sect:comment}
%
%This interesting question remains rather difficult to treat in generality.
%On the one hand, the fact that $\rho$ is not integrable would
%tend to make the solutions go to zero asymptotically, see Theorem \ref{thm:main.isot.0} in the case of bounded initial data.
%But on the other hand, since $u_0$ is unbounded, the solution also tends to go to infinity
%(see explicit examples in \cite{BrandleChasseigneFerreira09} for $\rho=1$).
%
%Thus, there is a balance between $\rho$, $J$, and the initial data $u_0$ which is not easy to handle and the question remains open.

\end{document}